\documentclass{ifacconf}

\usepackage{graphicx}      
\usepackage{natbib}        

\def\R{{\ifmmode{\rm I}\mkern-3.5mu{\rm R}
    \else\leavevmode\hbox{I}\kern-.16em\hbox{R}\fi}}
\def\Z{{\ifmmode{\rm Z}\mkern-7.2mu{\rm Z}
    \else\leavevmode\hbox{Z}\kern-.38em \hbox{Z}\fi }}
\def\N{{\ifmmode{\rm I}\mkern-3.5mu{\rm N}
    \else\leavevmode\hbox{I}\kern-.16em \hbox{N}\fi} }
\def\C{\ifmmode{{\rm C}\mkern-15mu{\phantom{\rm t}\vrule}}\mkern10mu
    \else\leavevmode\hbox{C}\kern-.5em\hbox{I}\kern.3em\fi}
\def\Q{\ifmmode{{\rm Q}\mkern-16mu{\phantom{\rm t}\vrule}}\mkern10mu
    \else\leavevmode\hbox{Q}\kern-.57em\hbox{I}\kern.3em\fi}

\begin{document}
\begin{frontmatter}

\title{Numerical experiments with multistep model-predictive control approaches 
  and sensitivity updates for the tracking control of cars\thanksref{footnoteinfo}} 

\thanks[footnoteinfo]{This material is based upon work supported by the Air Force Office of 
  Scientific Research, Air Force Materiel Command, USAF, under Award No, FA9550-14-11-0298.} 

\author[First]{Matthias Gerdts} 

\address[First]{Institute of Mathematics and Applied Computing, 
  Department of Aerospace Engineering, University of the Federal Armed Forces at Munich, 
  85577 Neubiberg, Germany (e-mail: matthias.gerdts@unibw.de)}

\begin{abstract}                
  The paper discusses multistep nonlinear model-predictive control (NMPC) schemes for the tracking 
  of a car model along a given reference track. In particular we will compare the numerical 
  performance and robustness of classic single step NMPC, multistep NMPC without re-optimization, 
  multistep NMPC with re-optimization, and multistep NMPC with sensitivity updates instead of 
  a full re-optimization. 
\end{abstract}

\begin{keyword}
  model predictive control, multistep NMPC, sensitivity updates, car tracking
\end{keyword}

\end{frontmatter}

\section{Introduction}


Nonlinear model predictive control (NMPC) is a feedback control paradigm with the capability to 
take into account control and/or state constraints, compare \cite{RawlingsMayne2009,Gru11} for a 
comprehensive overview and analysis. As such it is very powerful, but it 
relies on the repeated solution of nonlinear optimization problems on a moving time horizon. 
Especially in online computations the solution of the latter often turns out to be the computational 
bottleneck of NMPC and efficient numerical techniques are required, compare 
\cite{DiehlBockSchloeder2005}. Still, it is often not possible to fully solve these optimization 
problems within a given time frame. If this is the case, several modifications 
of the classic NMPC scheme exist. For instance, one could simply stop the iterative 
optimization procedure as soon as a time budget is consumed and accept the 
so far obtained result. Another way is to use {\em multistep NMPC schemes}, which do not just 
implement one control step of a computed solution, but more than one steps. 
This approach gains additional time to solve an optimization problem on a predicted 
preview horizon. On the downside, since deviations from the reference trajectory 
are not corrected at these steps, the approach is less robust than the 
classic NMPC scheme. To overcome this drawback, a re-optimization on the remaining part of 
the preview horizon can be performed in order to react on intermediate deviations. 
This leads to a {\em multistep NMPC scheme with re-optimization}. Finally, instead of performing 
a re-optimization, one could instead use parametric sensitivity analysis to update
the optimal solution on the preview horizon in the presence of perturbations. 
This leads to a {\em multistep NMPC scheme with sensitivity updates}, compare \cite{Zavala2008}. 
The purpose of the paper is to investigate and compare the classic scheme and the three modifications 
of the classic NMPC scheme in view of their tracking performance and numerical robustness. 
A theoretical investigation can be found in \cite{Palma2015}. 

Throughout, the aim is to construct a feedback control law $\mu:\N\times X \longrightarrow U$ for the constrained 
control system in discrete time 
\begin{eqnarray*}
  x(k+1) & = & f(x(k),u(k)), \qquad k=0,1,2,\ldots, \\
  x(k)   & \in & X, \qquad k=0,1,2,\ldots, \\
  u(k)   & \in & U, \qquad k=0,1,2,\ldots, \\
  x(0)   & = & x_0,
\end{eqnarray*} 
to track a given reference trajectory $(x_r(k),u_r(k))$, $k=0,1,2,\ldots$. 
Herein, $X\subset \R^n$ and $U\subset \R^m$ are given sets.
 
Often, the control system in discrete time can be interpreted as a discretization of a continuous process,
where $(x(k),u(k))$ corresponds to the state and the control at time $t_k = k h$
with sampling time $h>0$. Each of the different NMPC schemes yields a feedback control law $\mu_{N,M}$, 
where $N$ denotes the {\em preview horizon} and $M$ the {\em control horizon}. 
Closing the loop by setting $u(k) = \mu_{N,M}(k,x(K_M(k)))$ with $K_M(k) \leq k$ yields the closed-loop system
\begin{eqnarray*}
  x(k+1) & = & f(x(k),\mu_{N,M}(k,x(K_M(k))), \quad k=0,1,2,\ldots, \\
  x(0)   & = & x_0.
\end{eqnarray*} 
The feedback control laws $\mu_{N,M}$ will be defined in Sections~\ref{Sec:NMPC} and \ref{Sec:Sensitivity} 
for the following NMPC versions: 
\begin{itemize}
\item
  standard one-step NMPC,
\item
  multistep NMPC,
\item
  multistep NMPC with re-optimization,
\item
  multistep NMPC with sensitivity updates.
\end{itemize}
Numerical experiments with these schemes are presented in Section~\ref{Sec:Experiments}.

\section{NMPC Schemes} \label{Sec:NMPC}

Each of the NMPC schemes require to solve optimal control problems
in discrete time on some time horizon $[k_0,k_0+N]$ of the following type:

{\em 
  OCP($k_0,x_0,N$): \quad Minimize
  \begin{displaymath}
    \sum\limits_{k=k_0}^{k_0+N-1} f_0(k,x(k),u(k)) 
  \end{displaymath}
  subject to the constraints
  \begin{eqnarray*}
    x(k+1) & =    & f(x(k),u(k)), \qquad k=k_0,\ldots,k_0+N-1,\\
    x(k)   & \in  & X, \qquad k=k_0,\ldots,k_0+N,\\
    u(k)   & \in  & U, \qquad k=k_0,\ldots,k_0+N-1,\\ 
    x(k_0)   & =    & x_0.
\end{eqnarray*}}

Herein, $k_0$ denotes the current time of the process, $x_0$ the current (measured or predicted) state, 
and the number $N$ is called preview horizon. Throughout we consider tracking type objective functions. 
To this end, let a reference trajectory $(x_r(k),u_r(k))$, $k=0,1,2,\ldots$, 
be given. The function $f_0$ is then defined by
\begin{displaymath}
  f_0(k,x,u) = \| x_r(k) - x \|_V^2 + \| u_r(k) - u \|_W^2 
\end{displaymath}
with weighted norms $\| y \|_V = \sqrt{ y^\top V y }$, $\| z \|_W = \sqrt{ z^\top W z }$, where 
$V$ and $W$ are symmetric and positive semi-definite matrices. Throughout it is assumed that 
OCP($k_0,x_0,N$) for any choice of $(k_0,x_0,N)$ is feasible and possesses an optimal solution 
$(\hat x(k),\hat u(k))$, $k=k_0,\ldots,k_0+N-1$ (for notational simplicity we omit $\hat x(k_0+N)$ throughout), 
which can be computed by standard techniques, compare~\cite{Gerdts2012}.
The problem of infeasibility could be addressed in practice by relaxation of constraints or by choosing $N$ 
sufficiently large. 

The classic NMPC algorithm reads as follows and it yields a feedback law $\mu_N=\mu_{N,1}:\N\times X \longrightarrow U$.

\begin{alg}[classic NMPC]\label{Alg:Classic}$\ $\\[-12pt]
\begin{itemize}
\item[(0)]
  Input: preview horizon $N$, reference trajectory \linebreak $(x_r(\cdot),u_r(\cdot))$, weight matrices $V$ and $W$. 
  Set $k=0$. 
\item[(1)]
  Measure state $x(k) \in X$ at time $k$.
\item[(2)]
  Solve OCP($k,x(k),N$) on time horizon $[k,k+N]$. Let $\hat u(k),\ldots,\hat u(k+N-1)$ be the optimal control. 
\item[(3)]
  Define the feedback control $\mu_N(k,x(k)) := \hat u(k)$ and apply it: 
  \begin{displaymath}
    x(k+1) = f(x(k),\mu_N(k,x(k)))
  \end{displaymath}
\item[(4)]
  Set $k \leftarrow k+1$ and go to (1).
\end{itemize}
\end{alg}

\medskip
\begin{rem}
  Note that the implementation of $\mu_N(k,x(k))$ typically 
  is delayed by some $\delta>0$, where $\delta$ denotes the time to solve OCP($k,x(k)$). 
  Alternatively, one could use the predicted state $x(k+1)$ in step (3) to 
  solve the next problem OCP($k+1,x(k+1),N$) already during the step from $k$ to $k+1$. However, 
  the predicted state $x(k+1)$ deviates usually from the measured state at $k+1$ and hence an
  update of the computed solution might become necessary. This could be achieved by re-optimization or 
  by sensitivity updates as in Section~\ref{Sec:Sensitivity}. 
\end{rem}

The classic NMPC scheme requires to solve the optimal control problem at each 
time instance. If this is too time consuming, then the following multistep NMPC scheme 
is useful to reduce the number of optimal control problems to be solved. The idea is to 
apply not just the control $\hat u(k)$ in step (3) but to apply $M \leq N$ controls 
$\hat u(k),\hat u(k+1),\ldots,\hat u(k+M-1)$. The number $M$ is called {\em control horizon}. 

\begin{alg}[$M$-multistep NMPC]\label{Alg:Multi}$\ $\\[-12pt]
\begin{itemize}
\item[(0)]
  Input: preview horizon $N$, reference trajectory \linebreak $(x_r(\cdot),u_r(\cdot))$, weight matrices $V$ and $W$, 
  control horizon $M\leq N$. Set $k=0$. 
\item[(1)]
  Measure state $x(k)\in X$ at time $k$.
\item[(2)]
  Solve OCP($k,x(k),N$) on time horizon $[k,k+N]$. Let $\hat u(k),\ldots,\hat u(k+N-1)$ be the optimal control. 
\item[(3)]
  Define the feedback control 
  \begin{displaymath}
    \mu_{N,M}(k+j,x(k)) := \hat u(k+j), \quad j=0,\ldots,M-1,
  \end{displaymath}
  and apply it for $j=0,\ldots,M-1$: 
  \begin{displaymath}
    x(k+j+1) = f(x(k+j),\mu_{N,M}(k+j,x(k))).
  \end{displaymath}
\item[(4)]
  Set $k \leftarrow k+M$ and go to (1).
\end{itemize}
\end{alg}

Note that the $M$-multistep NMPC scheme and the classic NMPC scheme coincide for $M=1$. 
Apparently, the number of optimal control problems to be solved is reduced by a factor of $1/M$ in the 
$M$-multistep method when compared to the classic scheme. This reduction leaves more time to 
solve the optimal control problems. On the other hand, the controlled system runs for a longer time in 
open-loop mode and thus may become less stable, since the actual state is measured only every 
$M$ steps. 

In order to account for the potential stability issues in Algorithm~\ref{Alg:Multi}, a re-optimization 
can be performed. In order to re-use the previously computed optimal solution on $[k,k+N]$ as an 
initial guess, the re-optimization will only be performed on the reduced time horizon 
$[k+j,N]$. 

\begin{alg}[$M$-multistep NMPC with re-optimization]\label{Alg:Reopt}$\ $\\[-24pt]
\begin{itemize}
\item[(0)]
  Input: preview horizon $N$, reference trajectory \linebreak $(x_r(\cdot),u_r(\cdot))$, weight matrices $V$ and $W$, 
  control horizon $M\leq N$. Set $k=0$. 
\item[(1)]
  For $j=0,\ldots,M-1$ do
  \begin{itemize}
  \item[(1a)]
    Measure state $x(k+j)\in X$ at time $k+j$.
  \item[(1b)]
    Solve OCP($k+j,x(k+j),N-j$) on time horizon $[k+j,k+N]$. Let $\hat u(k+j),\ldots,\hat u(k+N-1)$ be the optimal control. 
  \item[(1c)]
    Define the feedback control 
    \begin{displaymath}
      \mu_{N,M}(k+j,x(k+j)) := \hat u(k+j)
    \end{displaymath}
    and apply it 
    \begin{displaymath}
      \hspace*{-6pt}x(k+j+1) = f(x(k+j),\mu_{N,M}(k+j,x(k+j))).
    \end{displaymath}
  \end{itemize}
\item[(2)]
  Set $k \leftarrow k+M$ and go to (1).
\end{itemize}
\end{alg}

\medskip
\begin{rem}
  Note that a re-optimization in step (1b) is only necessary, if the measured state 
  at $k+j$ deviates from the optimal state $\hat x(k+j)$ of the problem OCP($k+j-1,x(k+j-1),N-j-1$).
\end{rem}

A modification of Algorithm~\ref{Alg:Reopt}, which avoids the solution of the optimal control problems 
in step (1b), is described in the following Section~\ref{Sec:Sensitivity}.

\section{Multistep NMPC with Sensitivity Updates} \label{Sec:Sensitivity} 

The idea of the multistep NMPC scheme with sensitivity updates is to 
avoid to solve OCP($k+j,x(k+j),N-j$) in step (1b) of Algorithm~\ref{Alg:Reopt}. 
Instead, the solution of OCP($k+j,x(k+j),N-j$) will be approximated by means of 
a so-called sensitivity update, which will be the result of a parametric sensitivity analysis
of the optimal control problems with respect to the initial states.

\subsection{Parametric Sensitivity Analysis}

In order to perform the parametric sensitivity analysis, it is convenient to 
view the optimal control problems as a parametric optimization 
problem of type 

{\em 
  NLP($p$): \quad Minimize
  \begin{displaymath}
    J(z,p) 
  \end{displaymath}
  with respect to $z\in \R^{n_z}$ subject to the constraints
  \begin{eqnarray*}
    H(z,p) & = & 0,\\
    G(z,p) & \leq & 0.
\end{eqnarray*}}

Herein, $p\in\R^{n_p}$ denotes a parameter, $J:\R^{n_z} \longrightarrow \R$, 
$H:\R^{n_z}\longrightarrow \R^{n_H}$, $G:\R^{n_z}\longrightarrow \R^{n_G}$ are at least twice continuously 
differentiable functions. We are interested in properties of the {\em solution mapping} or {\em parameter-to-solution mapping} 
$p\mapsto z^*(p)$, where $z^*(p)$ denotes an optimal solution of NLP($p$). Conditions 
under which the solution mapping $z^*$ depends in a continuously differentiable way on the 
parameter $p$ are of particular interest, since in this case a linearization 
\begin{displaymath}
  z^*(p) = z^*(\hat p) + (z^*)'(\hat p) ( p - \hat p) + o(\|p-\hat p\|)
\end{displaymath}
around a nominal parameter $\hat p$ becomes possible. Neglecting the error term yields the 
approximate optimal solution $\tilde z(p)$ for $p$ close to $\hat p$: 
\begin{displaymath}
  z^*(p) \approx \tilde z(p) := z^*(\hat p)  + (z^*)'(\hat p) ( p - \hat p ), 
\end{displaymath}
compare \cite{BueskensMaurer2001}. 
For $p$ sufficiently close to $\hat p$, $\tilde z(p)$ may serve as a sufficiently good 
approximation to the optimal solution $z^*(p)$ of the perturbed nonlinear optimization 
problem NLP($p$). Note that the evaluation of $\tilde z(p)$ requires only a matrix 
vector multiplication and two vector additions, that is, the computational effort for these
operations is negligible. It remains to establish the solution differentiability and 
the computation of the sensitivity matrix $(z^*)'(\hat p)$. 

The solution differentiability of the map $z^*$ was established by \cite{Fiacco1983} 
with the following sensitivity theorem. The index set of active inequality constraints is given by 
\begin{displaymath}
  A(z,p) := \{ i  \; |\; G_i(z,p)=0, i\in \{1,\ldots,n_G\} \}.
\end{displaymath}
A local minimum $\hat z$ of NLP($\hat p$) is called {\em strongly regular}, 
if the following properties hold:
\begin{itemize}
\item[(a)]
  $\hat z$ fulfills the linear independence constraint qualification (LICQ), i.e. 
  the gradients $\nabla_z G_i(\hat z,\hat p)$, $i\in A(\hat z,\hat p)$, and $\nabla_z H_j(\hat z,\hat p)$, 
  $j=1,\ldots,n_H$, are linearly independent. 
\item[(b)]
  The KKT conditions hold at $(\hat z,\hat \mu,\hat\lambda)$, i.e. 
  \begin{displaymath}
    0 = \nabla_z L(\hat z,\hat \mu,\hat \lambda,\hat p), \ \hat\mu \geq 0, \hat \mu^\top G(\hat z,\hat p) = 0,
  \end{displaymath}
  where 
  \begin{displaymath}
    L(z,\mu,\lambda,p) := J(z,p) + \mu^\top G(z,p) + \lambda^\top H(z,p)
  \end{displaymath}
  denotes the Lagrange function of NLP($p$) with Lagrange multipliers $\mu$ and $\lambda$. 
\item[(c)]
  The strict complementarity condition holds: 
  \begin{displaymath}
    \hat \mu_i - G_i(\hat z,\hat p) > 0 \quad \mbox{for all } i=1,\ldots,n_G.
  \end{displaymath}
\item[(d)]
  We have
  \begin{displaymath}
    d^\top \nabla^2_{zz}L(\hat z,\hat \mu,\hat\lambda,\hat p) d > 0
  \end{displaymath}
  for all $d \in T_C(\hat z,\hat p)$, $d\not=0$, where 
  \begin{displaymath}
    \hspace*{-5pt}T_C(z,p) = 
    \left\{ d\; \left|\; 
    \begin{array}{rcl}       
      \displaystyle \nabla_z G_i(z,p)^\top d & = & 0, \; 
      i\in A(z,p), \\  
      \displaystyle \nabla_z H_j(z,p)^\top d & = & 0,\; j=1,\ldots,n_H
    \end{array}\right.\right\}.
  \end{displaymath}
\end{itemize}

A proof of the following theorem can be found in \cite{Fiacco1983} or \cite[Theorem 6.1.4]{Gerdts2012}.

\begin{thm}\label{Thm:Sensitivity} 
  Let $J$, $G$, and $H$ be twice continuously differentiable and $\hat p$ a fixed nominal 
  parameter. Let $\hat z$ be a strongly regular local minimum of NLP($\hat p$) with 
  Lagrange multipliers $\hat \lambda$ and $\hat \mu$. 

  Then there exist neighborhoods $B_{\epsilon}(\hat p)$ and  $B_{\delta}(\hat z,\hat \mu,\hat \lambda)$, 
  such that NLP($p$) has a unique strongly regular local minimum 
  \begin{displaymath}
    (z^*(p),\mu^*(p),\lambda^*(p)) \in B_{\delta}(\hat z,\hat \mu,\hat \lambda)
  \end{displaymath} 
  for each $p\in B_{\epsilon}(\hat p)$, and $A(\hat z,\hat p)=A(z^*(p),p)$.
 
  In addition, $(z^*(p),\mu^*(p),\lambda^*(p))$ is continuously 
  differentiable with respect to $p$ in these neighborhoods with 
  \begin{equation} \label{EQ_SENS_1}
    \left(
    \begin{array}{ccc}
      \nabla^2_{zz} L & (G_z')^\top & (H_z')^\top
      \\[6pt]
      \displaystyle
      \hat{\Xi}\cdot G_z' & \displaystyle \hat{\Gamma} & 0 \\[6pt] 
      \displaystyle
      H_z' & 0 & 0 
    \end{array}
  \right)\hspace*{-5pt}
  \begin{array}{rcl}
    \left(
  \begin{array}{c}\displaystyle
    (z^*)'(\hat{p}) \\[12pt] \displaystyle
    (\mu^*)'(\hat p) \\[12pt] \displaystyle
    (\lambda^*)'(\hat{p})
  \end{array}
  \right) =  -\left(\begin{array}{c}\displaystyle 
    \nabla^2_{zp} L \\[12pt] 
    \displaystyle \hat{\Xi} \cdot G_p'  \\[12pt] 
    \displaystyle H_p' 
  \end{array}\right)
  \end{array}
  \end{equation}
  where $\hat{\Xi} =  \mbox{{\rm diag}}(\hat{\mu}_1,\ldots,\hat{\mu}_{n_G})$, 
  $\hat{\Gamma} =  \mbox{{\rm diag}}(G_1,\ldots,G_{n_G})$. Herein, all 
  functions and their derivatives are evaluated at 
  $(\hat z,\hat \mu,\hat \lambda,\hat p)$. 
\end{thm}

\begin{rem}
  Please note, that all the assumptions needed to establish solution 
  differentiability can be checked numerically. Results without the strict 
  complementarity condition are derived in \cite{Jittorntrum1984}. 
\end{rem}

\subsection{Application to OCP in $M$-multistep NMPC}

For the application of Theorem~\ref{Thm:Sensitivity} 
we assume that the sets $X$ and $U$ in the problems OCP($k,x,N$) are defined by finitely many inequalities. 
We exploit the parametric sensitivity analysis to avoid the solution of OCP($k+j,x(k+j),N-j$) 
in step (1b) of Algorithm~\ref{Alg:Reopt}. Instead we approximate its solution 
by a sensitivity update and arrive at the following algorithm: 

\pagebreak
\begin{alg}[$M$-multistep NMPC with sensitivity upd.]\label{Alg:Sens}$\ $\\[-24pt]
\begin{itemize}
\item[(0)]
  Input: preview horizon $N$, reference trajectory \linebreak $(x_r(\cdot),u_r(\cdot))$, 
  weight matrices $V$ and $W$, control horizon $M\leq N$. Set $k=0$. 
\item[(1)]
  Measure state $x(k)\in X$ at time $k$ and solve OCP($k,x(k),N$) on time horizon $[k,k+N]$.
  Let $(\hat x(k+\ell),\hat u(k+\ell))$, $\ell=0,\ldots,N-1$, denote the optimal solution.
\item[(2)]
  Perform in parallel: 
  \begin{itemize}
  \item[(2a)]
    Define the feedback control $\mu_{N,M}(k,x(k)) := \hat u(k)$ and apply it 
    \begin{displaymath}
      x(k+1) = f(x(k),\mu_{N,M}(k,x(k))).
    \end{displaymath}
  \item[(2b)]
    For each $j=1,\ldots,M$ perform a sensitivity analysis of OCP($k+j,\hat x(k+j),N-j$) with respect to 
    the parameter $\hat p_j := \hat x(k+j)$. 

    Let $u^*_j(k+\ell)(\cdot)$, $\ell=j,\ldots,N-1$, denote the 
    solution mappings according to Theorem~\ref{Thm:Sensitivity}. 

    Let 
    \begin{displaymath}
      S_j := u_j^*(k+j)'(\hat p_j)
    \end{displaymath} 
    denote the sensitivity of the nominal control $\hat u(k+j)$ of OCP($k+j,\hat x(k+j),N-j$) 
    with respect to $p_j$ at $\hat p_j$. 
  \end{itemize}
\item[(3)]
  For $j=1,\ldots,M-1$ do
  \begin{itemize}
  \item[(3a)]
    Measure state $x(k+j)\in X$ at time $k+j$.
  \item[(3b)]
    Define the feedback control 
    \begin{eqnarray*}
      && \mu_{N,M}(k+j,x(k+j)) \\
      && \quad := \hat u(k+j) + S_j \cdot \left( x(k+j) - \hat x(k+j) \right)
    \end{eqnarray*}
    and apply it 
    \begin{displaymath}
      x(k+j+1) = f(x(k+j),\mu_{N,M}(k+j,x(k+j))).
    \end{displaymath}
  \end{itemize}
\item[(4)]
  Set $k \leftarrow k+M$ and go to (1).
\end{itemize}
\end{alg}

Some remarks are in order. Firstly, the sensitivity analysis is only justified under the assumptions 
of Theorem~\ref{Thm:Sensitivity} for sufficiently small perturbations $x(k+j)\approx \hat x(k+j)$. 
As a consequence the sensitivity analysis might not provide good approximations for large deviations. 
In the latter situation, it is recommended to fully re-solve OCP($k+j,\hat x(k+j),N-j$) as in Algorithm~\ref{Alg:Reopt}.
Still, the updated control in (3b) may serve as an initial guess. 

Secondly, it is important to note that the tails $(\hat x(k+\ell),\hat u(k+\ell))$, $\ell=j,\ldots,N-1$, 
of the optimal solution $(\hat x(k+\ell),\hat u(k+\ell))$, $\ell=0,\ldots,N-1$, of OCP($k,x(k),N$) 
are optimal for the problems OCP($k+j,\hat x(k+j),N-j$), $j=1,\ldots,M$, 
according to Bellman's optimality principle. Hence, by solving OCP($k,x(k),N$) in step (1), 
all nominal solutions to the problems OCP($k+j,\hat x(k+j),N-j$), $j=1,\ldots,M$, in step 
(2b) are known and the initial state of OCP($k+j,\hat x(k+j),N-j$), i.e. $\hat p_j=\hat x(k+j)$, 
can be viewed as a parameter entering the problem.   

Thirdly, please note that the parametric sensitivity analysis in step (2b) yields 
different solution mappings $u^*_j(k+\ell)(\cdot)$, $\ell=j,\ldots,N-1$, for each $j$. 
Only the first sensitivity $u^*_j(k+j)'(\cdot)$ at the time point $k+j$ is actually exploited in (3b).

It remains to compute the sensitivities $S_j$, $j=1,\ldots,M$, in step (2b) in an efficient way. 
The straight forward way of doing this is to solve equation (\ref{EQ_SENS_1}) for 
each of the problems OCP($k+j,\hat x(k+j),N-j$), $j=1,\ldots,M$, using the nominal solution
$(\hat x(k+\ell),\hat u(k+\ell))$, $\ell=j,\ldots,N-1$, and the nominal parameter $\hat p_j = \hat x(k+j)$. 
Please note that the dimension of the linear equation shrinks with increasing $j$ since 
the variables $(\hat x(k+\ell),\hat u(k+\ell))$, $\ell=0,\ldots,j-1$, and the constraints 
at the time points $k+\ell$, $\ell=0,\ldots,j-1$, are not present in OCP($k+j,\hat x(k+j),N-j$).

Since only the first sensitivity $u^*_j(k+j)'(\cdot)$ is actually exploited in (3b), solving the full 
linear systems is not necessary and an alternative and more efficient way is outlined in the sequel. 
Herein, the sensitivity analysis is merely performed for OCP($k,\hat x(k),N$) with respect to 
the parameter $\hat p_0 = \hat x(k)$. This yields the sensitivity differentials 
\begin{equation}\label{EQ:2}
  u_0^*(k+\ell)'(\hat p_0),\quad \ell=0,\ldots,N-1,
\end{equation}
at the time points $k+\ell$, $\ell=0,\ldots,N-1$, by solving equation (\ref{EQ_SENS_1}) once 
for the nominal solution $(\hat x(k+\ell),\hat u(k+\ell))$, $\ell=0,\ldots,N-1$. 
If a deviation $p_0=x(k)$ of $\hat p_0$ is detected, the optimal control can be updated by Taylor 
approximation 
\begin{displaymath}
  u_0^*(k+\ell)(p_0) \approx \hat u(k+\ell) + u_0^*(k+\ell)'(\hat p_0) ( p_0 - \hat p_0)
\end{displaymath}
for $\ell=0,\ldots,N-1$ and the state can be predicted by 
\begin{displaymath}
  x(k+\ell+1) = f(x(k+\ell),u_0^*(k+\ell)(p_0))
\end{displaymath}
for $\ell=0,\ldots,N-1$. 

Unfortunately, the $M$-multistep NMPC algorithm with sensitivity updates requires the 
sensitivities $S_j = u_j^*(k+j)'(\hat p_j)$ with $\hat p_j = \hat x(k+j)$ for $j=1,\ldots,M$, 
and not the sensitivities $u^*_0(k+j)'(\hat p_0)$ in (\ref{EQ:2}). Hence, a way to compute the $S_j$'s 
from (\ref{EQ:2}) is sought. To this end, we exploit the dynamics 
\begin{equation}\label{EQ:3}
  \hat x(k+1) = f(\hat x(k),\hat u(k)).
\end{equation}

\begin{assum}\label{Ass:1}
  Let the Jacobian matrix $f'_x(\hat x(k),\hat u(k))$ be non-singular.
\end{assum}

\begin{rem}
  Please note that Assumption~\ref{Ass:1} is satisfied for a sufficiently small step-size $h$, if 
  the dynamics are given by a one-step discretization method (e.g. a Runge-Kutta method) 
  for a differential equation, i.e. if $f$ is of type $f(x,u) = x + h \Phi(x,u,h)$. 
\end{rem}

If Assumption~\ref{Ass:1} holds, Equation (\ref{EQ:3}) can be solved for $\hat x(k)$ by the 
implicit function theorem, which yields the existence of 
neighborhoods $B_\epsilon(\hat x(k+1))$ and $B_\delta(\hat x(k))$ with $\epsilon>0$, $\delta >0$ 
and a mapping 
\begin{displaymath}
  \xi_0 : B_\epsilon(\hat x(k+1)) \longrightarrow B_\delta(\hat x(k))
\end{displaymath}
such that $\hat x(k) = \xi_0(\hat x(k+1))$ and 
\begin{displaymath}
  x(k+1) = f(\xi_0(x(k+1)),\hat u(k))
\end{displaymath}
holds for every $x(k+1)\in B_\epsilon(\hat x(k+1))$. 
Moreover, by differentiating this identity with respect to $x(k+1)$ we find 
\begin{eqnarray*}
  I & = & f'_x(\xi_0(\hat x(k+1)),\hat u(k)) \cdot \xi_0'(\hat x(k+1))
\end{eqnarray*}
and thus 
\begin{eqnarray*}
  \xi_0'(\hat x(k+1)) & = & f'_x(\xi_0(\hat x(k+1)),\hat u(k))^{-1} \\
  & = & f'_x(\hat x(k),\hat u(k))^{-1}.
\end{eqnarray*}
Note that $\xi_0'(\hat x(k+1))$ is the derivative of the initial state $\hat x(k)$ with 
respect to $x(k+1)$. 

Note further, that we have the relation 
\begin{displaymath}
  u^*_1(k+1) ( x(k+1)) = u^*_0(k+1)( \xi_0(x(k+1)) )
\end{displaymath}
for every $x(k+1)\in B_\epsilon(\hat x(k+1))$ (eventually after reducing $\epsilon$ taking into 
account the neighborhoods of the sensitivity theorem).  
 
Now, by the chain rule we obtain
\begin{eqnarray*}
  S_1 & = & u_1^*(k+1)'(\hat x(k+1)) \\
  & = & u_0^*(k+1)'(\hat x(k)) \cdot \xi'_0(\hat x(k+1)) \\
  & = & u_0^*(k+1)'(\hat x(k)) \cdot f'_x(\hat x(k),\hat u(k))^{-1}. 
\end{eqnarray*}
This formula allows to compute $S_1$ without solving equation (\ref{EQ_SENS_1}) 
for OCP($k+1,\hat x(k+1),N-1$). 

This construction can be repeated for $j=2,\ldots,M$ exploiting the relations
\begin{eqnarray*}
  && u^*_j(k+j) ( x(k+j) ) \\
  && \quad = u^*_0(k+j)( \xi_0 \circ \xi_1 \cdots \circ \xi_{j-1}(x(k+j)) ), 
\end{eqnarray*}
where $\xi_{j-1}$ satisfies $\hat x(k+j-1) = \xi_{j-1}(\hat x(k+j))$ and 
\begin{displaymath}
  x(k+j) = f(\xi_{j-1}(x(k+j)),\hat u(k+j-1))
\end{displaymath}
holds for every $x(k+j)$ in some neighborhood of $\hat x(k+j)$. Herein, 
Assumption \ref{Ass:1} has to hold accordingly for $f'_x(\hat x(k+\ell),\hat u(k+\ell))$, 
$\ell = 1,\ldots,M-1$.  Then we obtain
\begin{eqnarray*}
  S_j & = & u_j^*(k+j)'(\hat x(k+j)) \\
  & = & u_0^*(k+j)'(\hat x(k)) \cdot \prod_{\ell=0}^{j-1} f'_x(\hat x(k+\ell),\hat u(k+\ell))^{-1}. 
\end{eqnarray*}

For a rigorous mathematical stability and performance analysis of the different $M$-multistep NMPC schemes 
we refer the reader to \cite{Palma2015}. 

\section{Numerical Experiments}\label{Sec:Experiments}

\begin{figure}
  \begin{center}
    \includegraphics[height=0.24\textwidth,angle=-90]{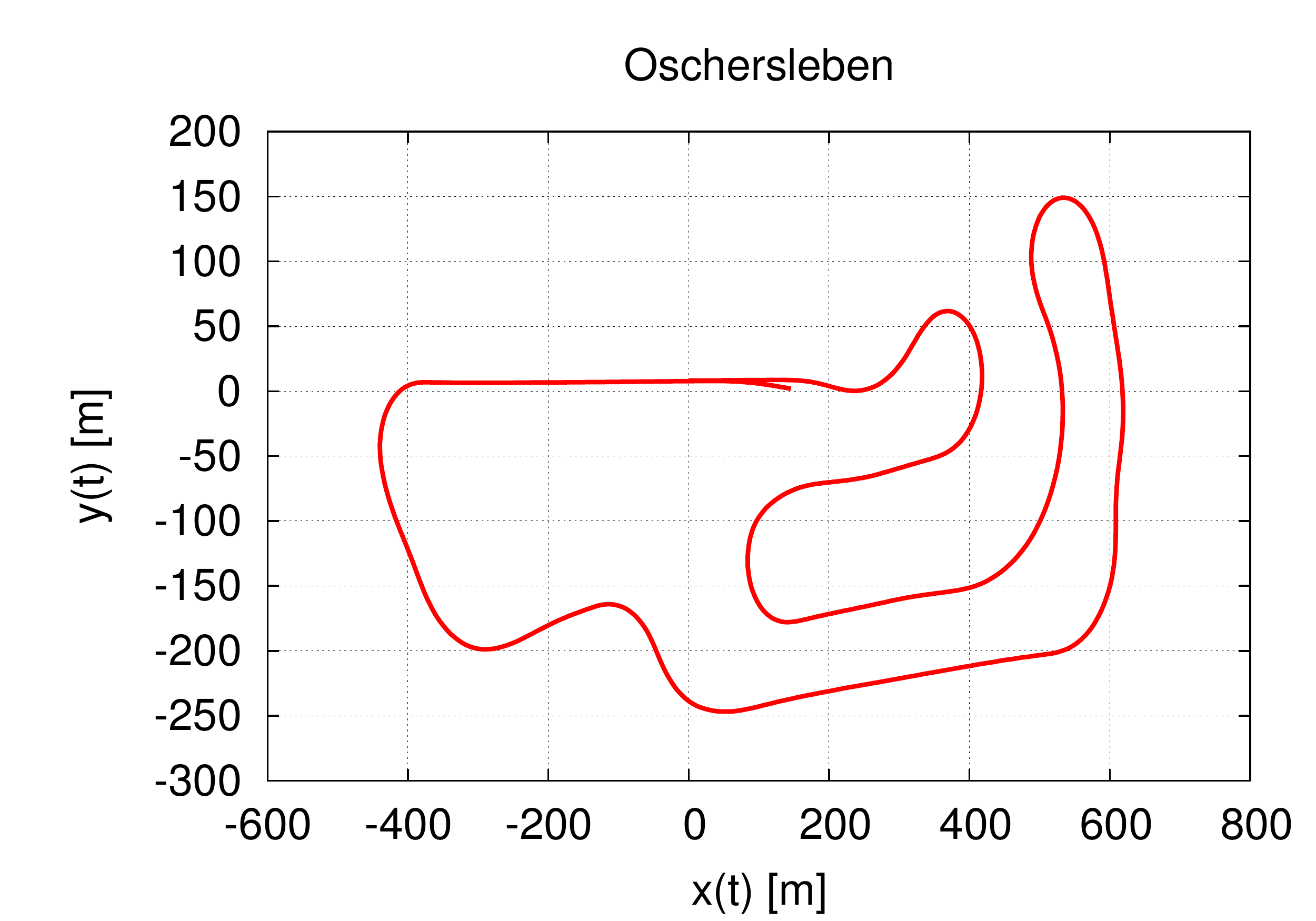}
    \includegraphics[height=0.24\textwidth,angle=-90]{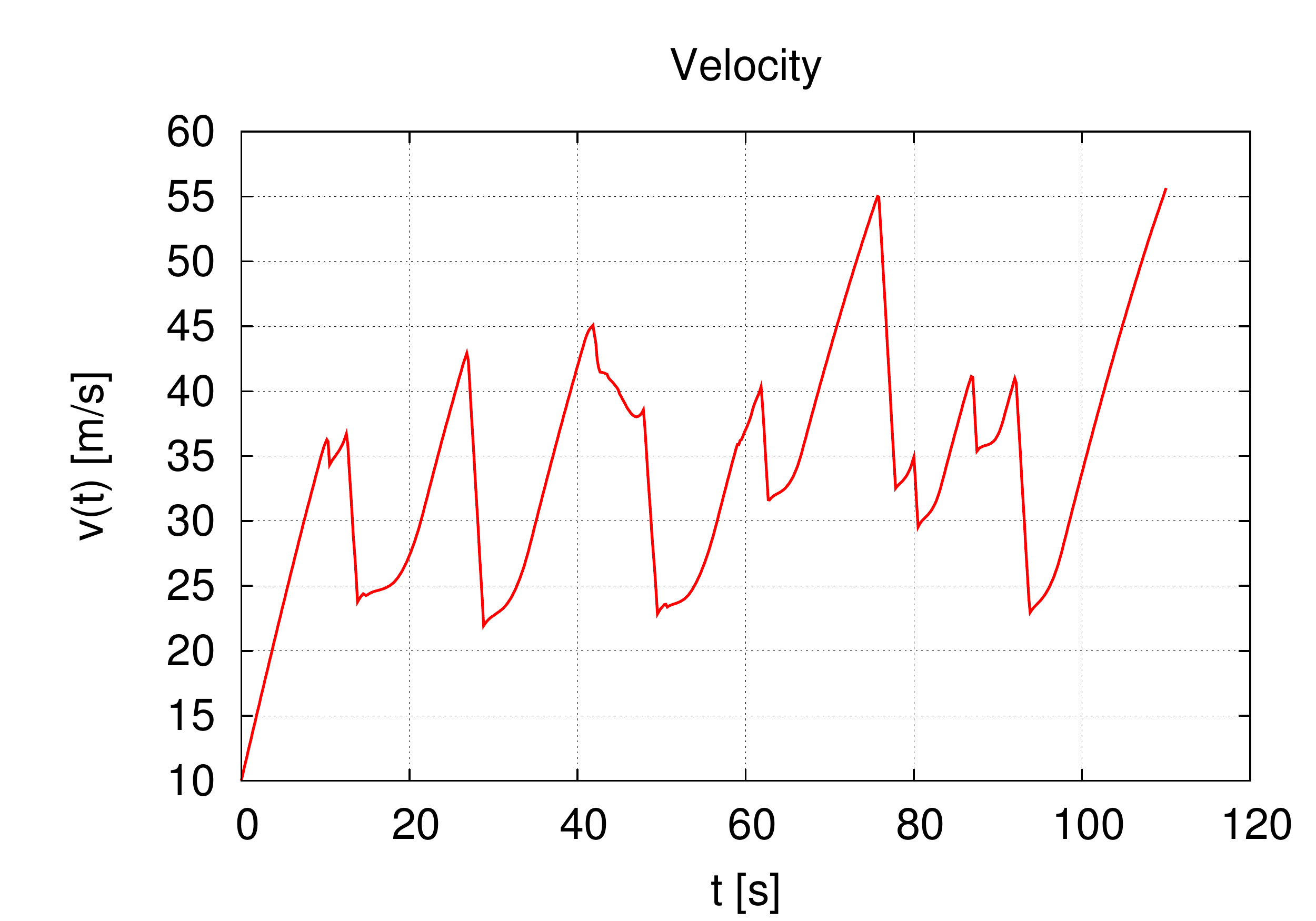}
    \includegraphics[height=0.24\textwidth,angle=-90]{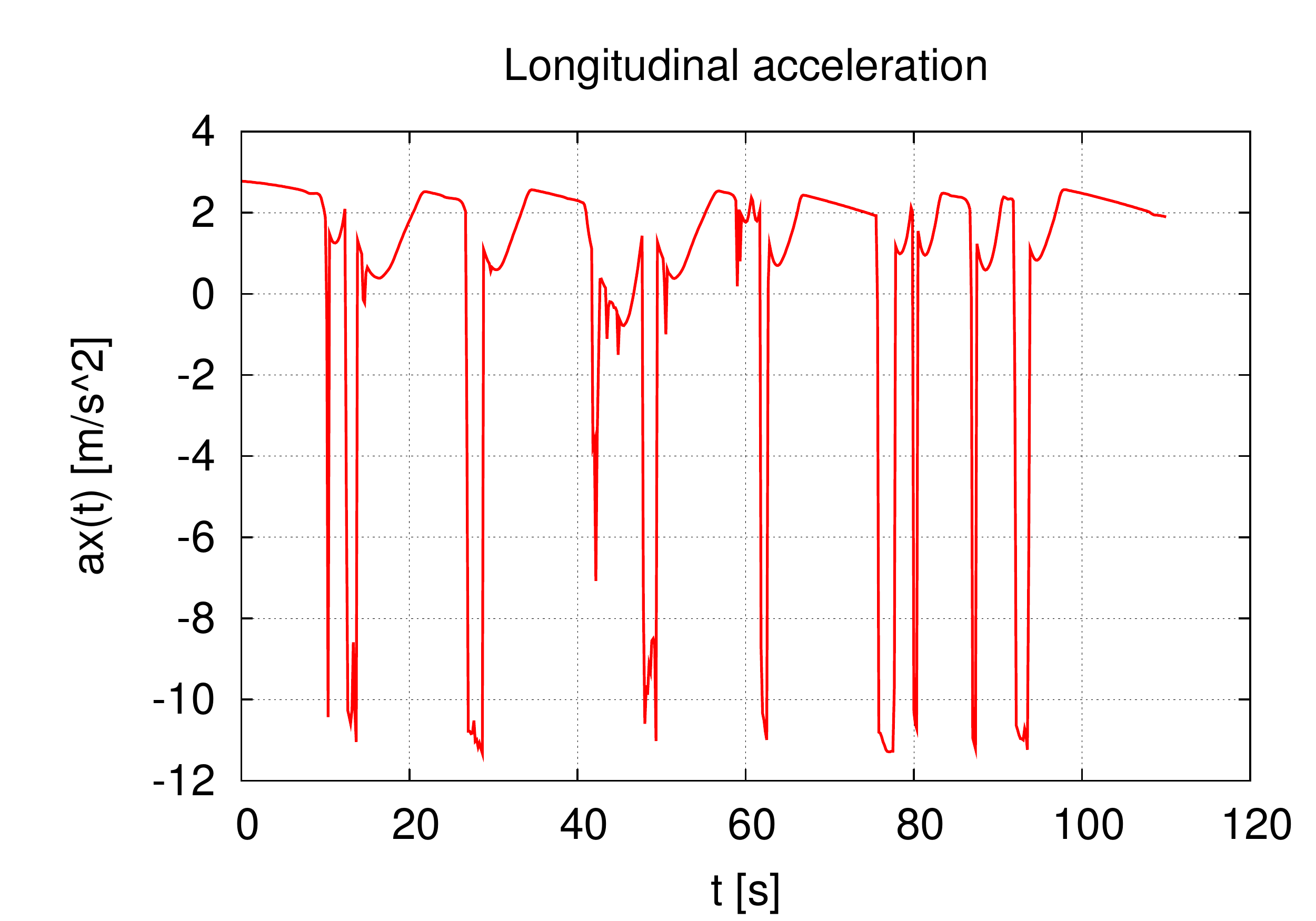}
    \includegraphics[height=0.24\textwidth,angle=-90]{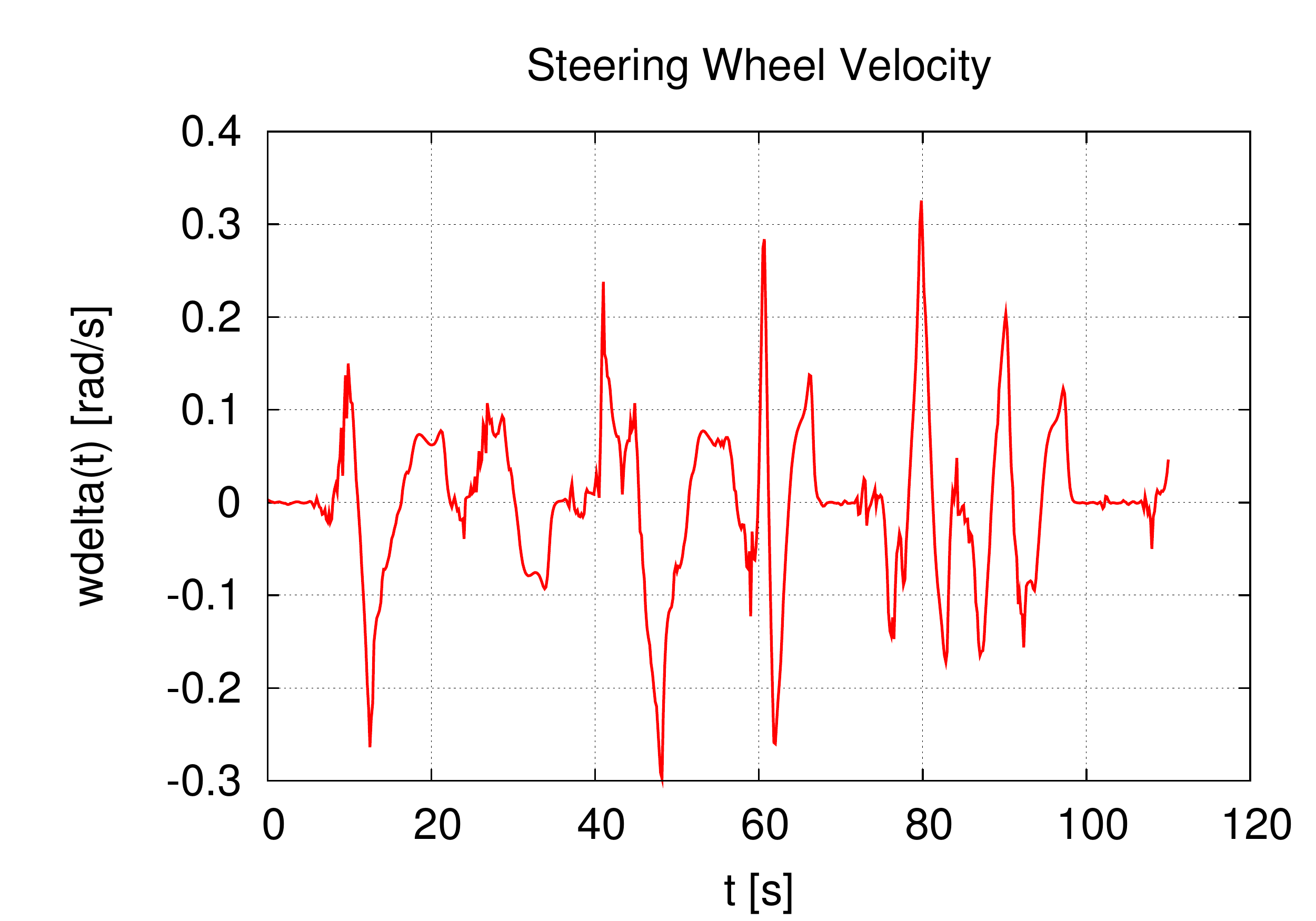}    
  \end{center}
  \caption{Reference trajectory for the racing track of Oschersleben: position $(x_r,y_r)$ (top left),
    velocity $v_r$ (top right), acceleration $u_{1,r}$ (bottom left), and steering angle velocity $u_{2,r}$ (bottom right).}\label{Fig:1}
\end{figure}

We compare the four NMPC schemes for the problem of tracking the raceline along the 
testtrack of Oschersleben in Figure~\ref{Fig:1} with the following kinematic car model: 
\begin{eqnarray*}
  x'(t) & = & v(t) \cos \psi(t), \qquad x(0) = x_0,\\
  y'(t) & = & v(t) \sin \psi(t), \qquad y(0) = y_0,\\
  \psi'(t) & = &\frac{v(t)}{\ell} \tan\delta(t), \qquad \psi(0) = \psi_0,\\
  v'(t) & = &u_1(t), \qquad v(0) = v_0,\\
  \delta'(t) & = &u_2(t), \qquad \delta(0) = \delta_0.
\end{eqnarray*}
Herein, $\ell = 4\ [m]$ denotes the length of the car, $(x,y)$ the position of the center 
of the rear axle, $\psi$ the yaw angle, $v$ the velocity, and $\delta$ the steering angle.
All numerical experiments have been conducted with $(x_0,y_0,\psi_0,v_0,\delta_0)=(0\ [m],0\ [m],0\ [rad],10\ [m/s],0\ [rad])$, 
a preview horizon of $T=3\ [s]$, $N=11$ grid points (i.e. a step-size of $h=0.3\ [s]$), 
control horizon $M=3$. The initial position 
$(x,y)$ and the velocity $v$ in each step of the NMPC schemes are perturbed by 
adding equally distributed noise in the range $[-0.05,0.05]$, which is realistic for measurements with 
a differential GPS system. The total control time horizon was $t_f=110\ [s]$. 
The controls are subject to the control bounds $u_1 \in [-12,3]\ [m/s^2]$ and 
$u_2\in [-0.5,0.5]\ [rad/s]$. 
Moreover, the state constraints $v \in [0,60]\ [m/s]$ and $\delta \in [-0.5,0.5]\ [rad]$ 
have to be obeyed. Throughout, the objective function 
\begin{eqnarray*}
 && \int_{0}^T \alpha_1 \left\| \left(\begin{array}{c}  x(t) -x_r(t) \\ y(t)-y_r(t) \end{array}\right) \right\|^2 + \alpha_2 (v(t) - v_r(t))^2 \\
 && \qquad + \alpha_3 \left\| \left(\begin{array}{c}  u_1(t) -u_{1,r}(t) \\ u_2(t)-u_{2,r}(t) \end{array}\right) \right\|^2 \ dt
\end{eqnarray*}
with $\alpha_1 = 1$, $\alpha_2 = 10^{-1}$, and $\alpha_3=10^{-3}$ was used in the NMPC schemes.  
The optimal control package {\ttfamily OCPID-DAE1} \cite{OCPIDDAE1} was used for solving the optimal control 
problems and performing the sensitivity analysis. The focus of the study is on the robustness and tracking error 
of the methods, not on the CPU times. For this reason, the sensitivities are computed by solving (\ref{EQ_SENS_1}) 
for simplicity. 

\begin{figure}
  \begin{center}
    \includegraphics[height=0.45\textwidth,angle=-90]{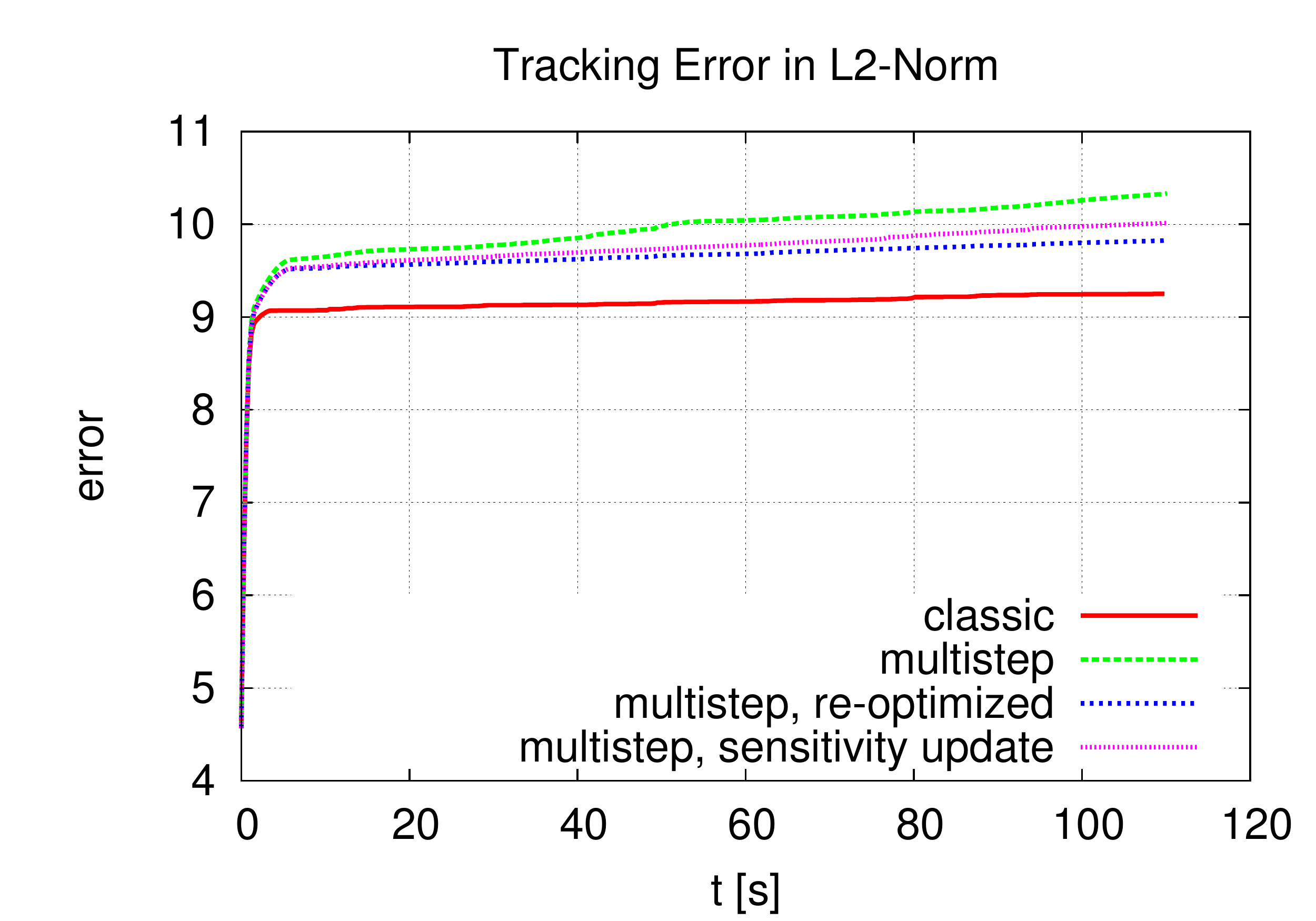}
  \end{center}
  \caption{Tracking error for the four different NMPC strategies.}\label{Fig:2}
\end{figure}

Figure~\ref{Fig:2} shows the tracking error $\| (x-x_r,y-y_r,v-v_r) \|_{L_2((0,t_f))}$ 
measured in the $L_2$-norm. The results show that the classic NMPC scheme performs best with 
regard to the tracking error, followed by the multistep NMPC scheme with re-optimization, 
the multistep NMPC scheme with sensitivity updates, and the multistep NMPC scheme. This outcome 
is the expected one since the classic scheme optimizes in each step on the full preview horizon
while the multistep scheme optimizes only after $M$ shifts have been performed. 
The multistep scheme with re-optimization re-optimizes at least on a shrinking horizon at 
every shift. The same holds true for the multistep scheme with sensitivity updates, but 
this only provides a Taylor approximation to the optimal solution. 
The large initial error is due to a large deviation of about $8.3\ [m]$ 
in the y-direction from the initial state of the reference solution

\begin{figure}
   \begin{center}
     \includegraphics[height=0.47\textwidth,angle=-90]{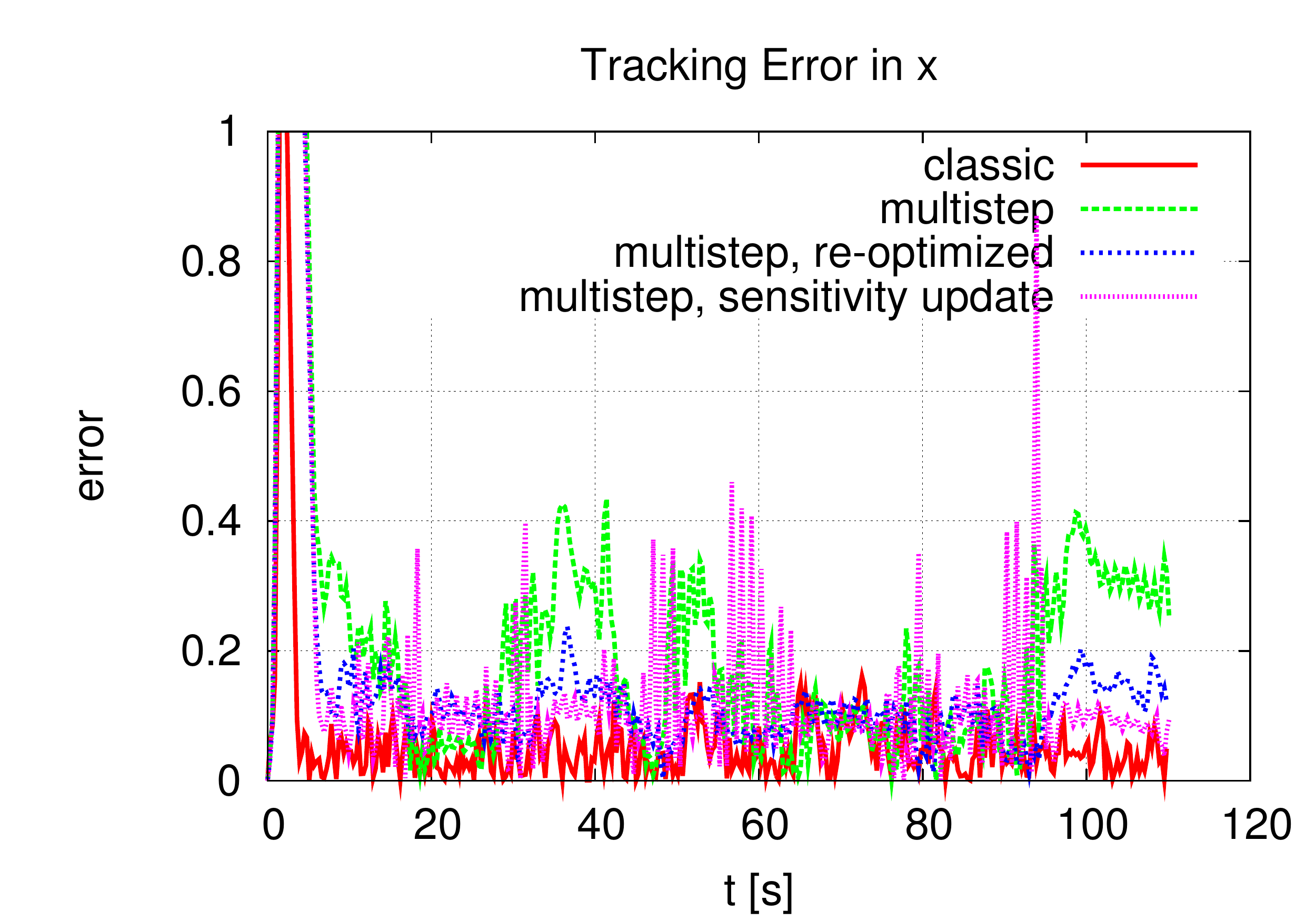}
     \includegraphics[height=0.47\textwidth,angle=-90]{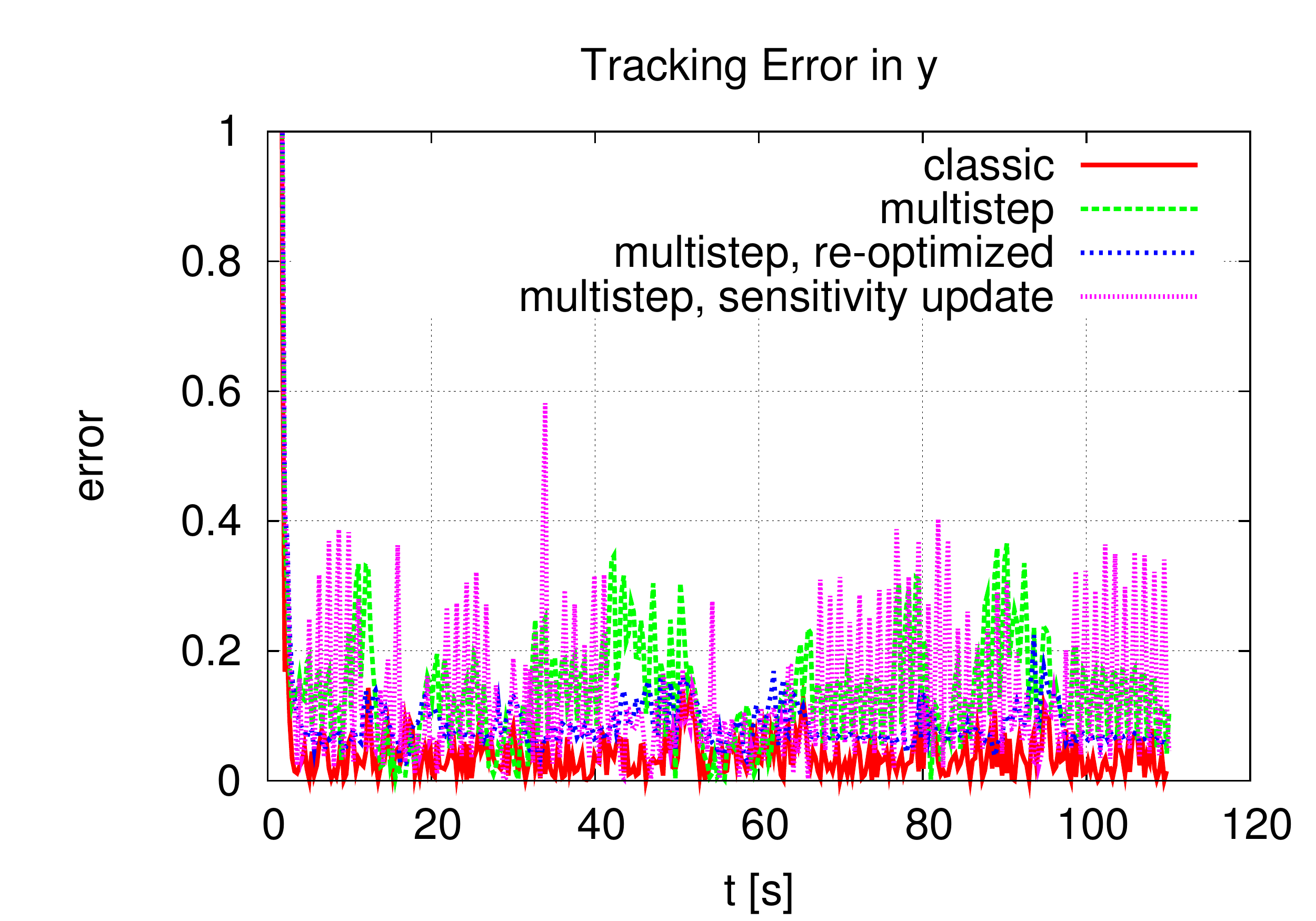} 
     \includegraphics[height=0.47\textwidth,angle=-90]{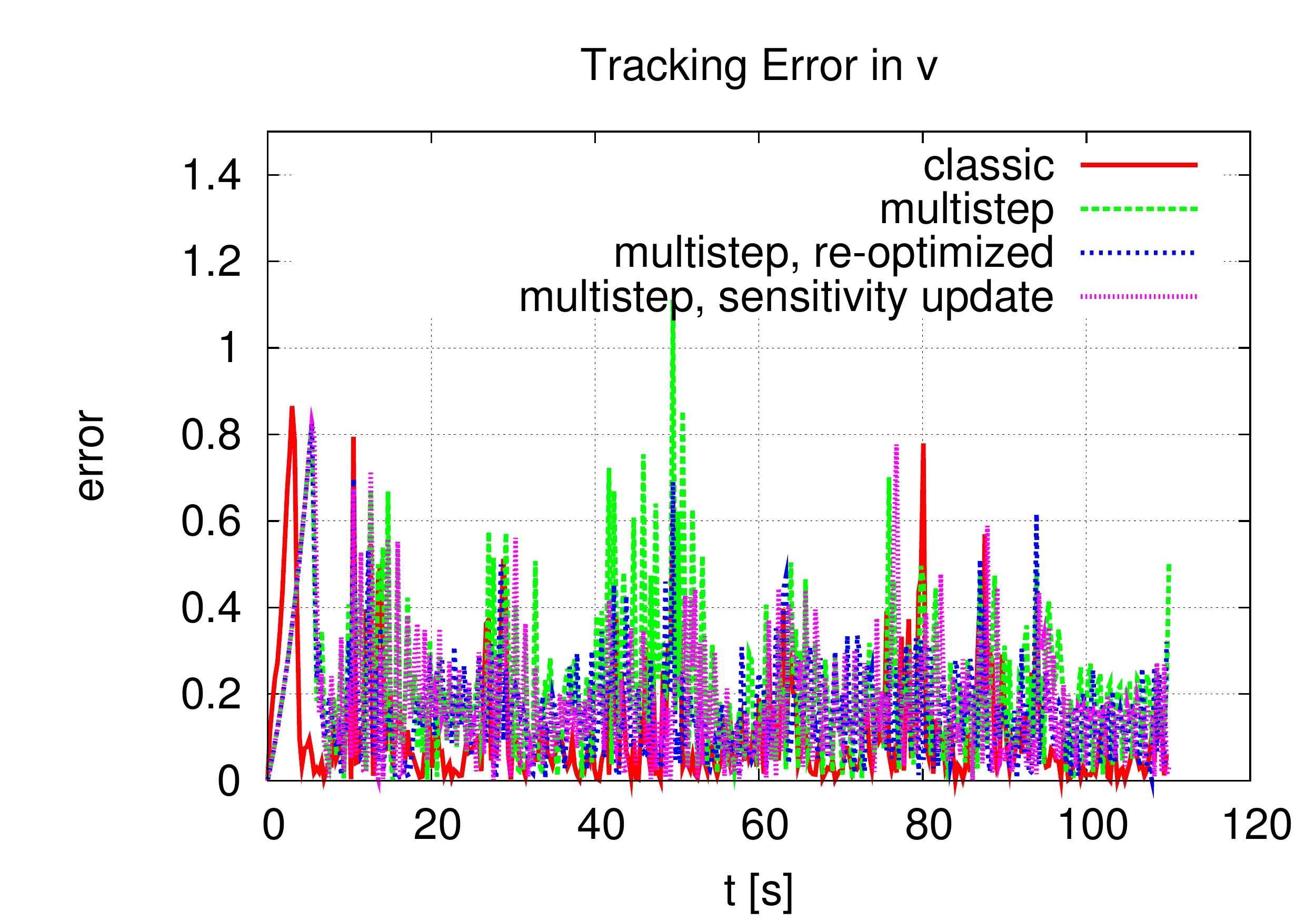}   
   \end{center}
   \caption{Comparison of the four NMPC schemes: error in the position x-position (top), y-position (middle), and velocity (bottom).}\label{Fig:3}
\end{figure}

Figure~\ref{Fig:3} shows the errors of the (x,y)-position and the velocity for the 
four NMPC schemes. All schemes are able to track the reference solution at a high precision.
Recall that equally distributed noise with an amplitude of 0.1 was added in the NMPC schemes.

\section{Conclusion}

A numerical study of the classic NMPC scheme, the multistep NMPC scheme, the multistep 
NMPC scheme with re-optimization, and the multistep NMPC scheme with sensitivity updates 
was performed and tested for a tracking problem along a racing track with a kinematic car 
model. The numerical study shows that all approaches are feasible and are able to track the 
given reference trajectory subject to random noise. Moreover, the results supports the 
expectation that the classic NMPC scheme performs best with regard to the tracking error. 
It is followed by the multistep scheme with re-optimization and the multistep scheme with 
sensitivity updates. Finally, the basic multistep scheme yields the largest tracking error of
the four approaches. 

\begin{ack}
  Many thanks to Johannes Goergen for the support with the implementation. 
\end{ack}

\bibliography{ifacconf}             

\end{document}